\documentstyle[11pt]{article}
\newlength{\standardunitlength}
\newtheorem{cor}{Corollary} \newtheorem{lemma}{Lemma}
\newtheorem{theorem}{Theorem}
\newenvironment{proof}{\noindent {\sc Proof:}}{$\Box$ \vspace{2 ex}}
\begin{document}

\begin{center} {\bf The Eigenvalue Spacing of a Random Unipotent
Matrix in its Representation on Lines} \end{center}

\begin{center}
By Jason Fulman
\end{center}

\begin{center}
Dartmouth College
\end{center}

\begin{center}
Hanover, NH 03755, USA
\end{center}

\begin{center}
Fulman@Dartmouth.Edu
\end{center}

\newpage
\begin{center}
Proposed Running Head: The Eigenvalue Spacing
\end{center}

\begin{abstract} The eigenvalue spacing of a uniformly chosen random finite
unipotent matrix in its permutation action on lines is studied. We
obtain bounds for the mean number of eigenvalues lying in a fixed arc
of the unit circle and offer an approach toward other asymptotics. For
the case of all unipotent matrices, the proof gives a probabilistic
interpretation to identities of Macdonald from symmetric function
theory. For the case of upper triangular matrices over a finite field,
connections between symmetric function theory and a probabilistic
growth algorithm of Borodin emerge. \end{abstract}

\begin{center} Key words: Eigenvalue spacing, symmetric functions,
Hall-Littlewood polynomial, random matrix.
\end{center}

\section{Introduction} The subject of eigenvalues of random matrices
is very rich. The eigenvalue spacings of a complex unitary matrix chosen from
Haar measure relate to the spacings between the zeros of the Riemann
zeta function (\cite{Odlyzko}, \cite{RS1}, \cite{RS2}). For further
recent work on random complex unitary matrices, see \cite{DS},
\cite{Rains}, \cite{Wieand}. The references \cite{Dyson} and
\cite{Mehta} contain much of interest concerning the eigenvalues of a
random matrix chosen from Dyson's orthogonal, unitary, and symplectic
circular ensembles, for instance connections with the statistics of
nuclear energy levels.

Little work seems to have been done on the eigenvalue spacings of matrices
chosen from finite groups. One recent step is Chapter 5 of Wieand's
thesis \cite{Wieand}. She studies the eigenvalue spacings of a random element
of the symmetric group in its permutation representation on the set
$\{1,\cdots,n\}$. This note gives two natural $q$-analogs
of Wieand's work. For the first $q$-analog, let $\alpha \in GL(n,q)$
be a random unipotent matrix. Letting $V$ be the vector space on which
$\alpha$ acts, we consider the eigenvalues of $\alpha$ in the
permutation representation of $GL(n,q)$ on the lines of $V$. Let
$X^{\theta}(\alpha)$ be the number of eigenvalues of $\alpha$ lying in
a fixed arc $(1,e^{i2\pi \theta}], 0 < \theta < 1$ of the unit
circle. Bounds are obtained for the mean of $X^{\theta}$ (we believe
that as $n \rightarrow \infty$ with $q$ fixed, a normal limit theorem
holds). A second $q$-analog which we analyze is the case when $\alpha$
is a randomly chosen unipotent upper triangular matrix over a finite field. A
third interesting $q$-analog would be taking $\alpha$ uniformly chosen
in $GL(n,q)$; however this seems intractable.

The main method of this paper is to interpret identities of symmetric
function theory in a probabilistic setting. Section \ref{ident} gives
background and results in this direction. This interaction appears
fruitful, and it is shown for instance that a probabilistic algorithm
of Borodin describing the Jordan form of a random unipotent upper
triangular matrix \cite{B} follows from the combinatorics of symmetric
functions. This ties in with work of the author on analogous
algorithms for the finite classical groups \cite{F}. The applications
to the eigenvalue problems described above appear in Section
\ref{applications}.  We remark that symmetric function theory plays
the central role in work of Diaconis and Shahshahani \cite{DS} on the
eigenvalues of random complex classical matrices.

\section{Symmetric functions} \label{ident}

	To begin we describe some notation, as on pages 2-5 of
\cite{Mac}. Let $\lambda$ be a partition of a non-negative integer $n
= \sum_i \lambda_i$ into non-negative integral parts $\lambda_1 \geq
\lambda_2 \geq \cdots \geq 0$. The notation $|\lambda|=n$ will mean
that $\lambda$ is a partition of $n$. Let $m_i(\lambda)$ be the number
of parts of $\lambda$ of size $i$, and let $\lambda'$ be the partition
dual to $\lambda$ in the sense that $\lambda_i' = m_i(\lambda) +
m_{i+1}(\lambda) + \cdots$. Let $n(\lambda)$ be the quantity $\sum_{i
\geq 1} (i-1) \lambda_i$. It is also useful to define the diagram
associated to $\lambda$ as the set of points $(i,j) \in Z^2$ such that
$1 \leq j \leq \lambda_i$. We use the convention that the row index
$i$ increases as one goes downward and the column index $j$ increases
as one goes across. So the diagram of the partition $(5441)$ is:

\[ \begin{array}{c c c c c}
		. & . & . & . & .  \\
		. & . & . & . &    \\
		. & . & . & . &    \\
		. & & & &  
	  \end{array}  \] Let $G_{\lambda}$ be an abelian $p$-group isomorphic to
$\bigoplus_i Cyc(p^{\lambda}_i)$. We write $G=\lambda$ if $G$ is an
abelian $p$-group isomorphic to $G_{\lambda}$. Finally, let
$(\frac{1}{p})_r = (1-\frac{1}{p}) \cdots (1-\frac{1}{p^r})$.

	The rest of the paper will treat the case $GL(n,p)$ with $p$
prime as opposed to $GL(n,q)$. This reduction is made only to make the
paper more accessible at places, allowing us to use the language of
abelian $p$-groups rather than modules over power series rings. From
Chapter 2 of Macdonald \cite{Mac} it is clear that everything works
for prime powers.

\subsection{Unipotent elements of $GL(n,p)$} \label{subunip} It is
well known that the unipotent conjugacy classes of $GL(n,p)$ are
parametrized by partitions $\lambda$ of $n$. A representative of the
class $\lambda$ is given by

\[ \left( \begin{array}{c c c c}
		M_{\lambda_1} & 0 & 0 &0 \\
		0 & M_{\lambda_2} & 0 & 0\\
		0 & 0 & M_{\lambda_3} & \cdots \\
		0 & 0 & 0 & \cdots
	  \end{array} \right), \] where $M_i$ is the $i*i$ matrix of the form

\[ \left( \begin{array}{c c c c c c}
		1 & 1 & 0 & \cdots & \cdots & 0\\
		0 & 1 & 1 & 0 & \cdots & 0\\
		0 & 0 & 1 & 1 & \cdots & 0\\
		\cdots & \cdots & \cdots & \cdots & \cdots & \cdots \\
		\cdots & \cdots & \cdots & 0 & 1 & 1\\
		0 & 0 & 0 & \cdots & 0 & 1
	  \end{array} \right). \]

 Lemmas \ref{classsize}-\ref{prod} recall elementary facts about
unipotent elements in $GL(n,p)$.

\begin{lemma} \label{classsize} (\cite{Mac} page 181,\cite{SS}) The
number of unipotent elements in $GL(n,p)$ with conjugacy class type
$\lambda$ is

\[ \frac{|GL(n,p)|}{p^{\sum (\lambda_i')^2} \prod_i
(\frac{1}{p})_{m_i(\lambda)}} .\] \end{lemma}

	Chapter 3 of \cite{Mac} defines Hall-Littlewood symmetric
functions $P_{\lambda}(x_1,x_2,\cdots;t)$ which will be used
extensively. There is an explicit formula for the Hall-Littlewood
polynomials. Let the permutation $w$ act on the $x$-variables by
sending $x_i$ to $x_{w(i)}$. There is also a coordinate-wise action of
$w$ on $\lambda=(\lambda_1, \cdots,\lambda_n)$ and $S^{\lambda}_n$ is
defined as the subgroup of $S_n$ stabilizing $\lambda$ in this
action. For a partition $\lambda=(\lambda_1,\cdots,\lambda_n)$ of
length $\leq n$, two formulas for the Hall-Littlewood polynomial
restricted to $n$ variables are:

\begin{eqnarray*}
P_{\lambda}(x_1,\cdots,x_n;t) & = & [\frac{1}{\prod_{i \geq 0}
\prod_{r=1}^{m_i(\lambda)} \frac{1-t^r}{1-t}}] \sum_{w \in S_n}
w(x_1^{\lambda_1} \cdots x_n^{\lambda_n} \prod_{i<j} \frac{x_i-tx_j}
{x_i-x_j})\\
& = & \sum_{w \in S_n/S_n^{\lambda}} w(x_1^{\lambda_1} \cdots x_n^{\lambda_n} \prod_{\lambda_i > \lambda_j} \frac{x_i-tx_j}{x_i-x_j})
\end{eqnarray*}

\begin{lemma} \label{formula} The probability that a unipotent element of $GL(n,p)$ has conjugacy class of type $\lambda$ is equal to either of

\begin{enumerate}
\item $\frac{p^n (\frac{1}{p})_n}{p^{\sum (\lambda_i')^2} \prod_i
(\frac{1}{p})_{m_i(\lambda)}}$

\item $\frac{p^n (\frac{1}{p})_n P_{\lambda}(\frac{1}{p},\frac{1}{p^2},\frac{1}{p^3},\cdots;\frac{1}{p})}{p^{n(\lambda)}}$
\end{enumerate}
\end{lemma}

\begin{proof} The first statement follows from Lemma \ref{classsize}
and Steinberg's theorem that $GL(n,p)$ has $p^{n(n-1)}$ unipotent
elements. The second statement follows from the first and from
elementary manipulations applied to Macdonald's principal
specialization formula (page 337 of \cite{Mac}). Full details appear
in \cite{F2}.  \end{proof}

	One consequence of Lemma \ref{formula} is that in the $p
\rightarrow \infty$ limit, all mass is placed on the partition
$\lambda=(n)$. Thus the asymptotics in this paper will focus on the
more interesting case of the fixed $p$, $n \rightarrow \infty$ limit.

\begin{lemma} \label{prod}

\[ \sum_{\lambda \vdash n} \frac{1}{p^{\sum (\lambda_i')^2} \prod_i
(\frac{1}{p})_{m_i(\lambda)}} = \frac{1}{p^n(\frac{1}{p})_n}\]

\end{lemma}

\begin{proof} Immediate from Lemma \ref{formula}.
\end{proof}

Lemmas \ref{duality} and \ref{likeMac} relate to the theory of Hall
polynomials and Hall-Littlewood symmetric functions \cite{Mac}. Lemma
\ref{duality}, for instance, is the duality property of Hall
polynomials.

\begin{lemma} \label{duality} (Page 181 of \cite{Mac}) For all
partitions $\lambda,\mu,\nu$,

\[|\{G_1 \subseteq
G_{\lambda}:G_{\lambda}/G_1=\mu,G_1=\nu\}| = |\{G_1 \subseteq
G_{\lambda}:G_{\lambda}/G_1=\nu,G_1=\mu\}| .\]
\end{lemma}

\begin{lemma} \label{likeMac} Let $G_{\lambda}$ denote an abelian $p$-group of type $\lambda$, and $G_1$ a subgroup. Then for all types $\mu$, 

\[ \sum_{\lambda \vdash n} \frac{\{|G_1 \subseteq
G_{\lambda}:G_1=\mu|\}}{p^{\sum (\lambda_i')^2} \prod_i
(\frac{1}{p})_{m_i(\lambda)}} = \frac{1}{p^{\sum (\mu_i')^2} \prod_i
(\frac{1}{p})_{m_i(\mu)}} \frac{1}{p^{n-|\mu|}
(\frac{1}{p})_{n-|\mu|}}.\] \end{lemma}

\begin{proof} Macdonald (page 220 of \cite{Mac}), using
Hall-Littlewood symmetric functions, establishes for any partitions
$\mu,\nu$, the equation:

\[ \sum_{\lambda: |\lambda|=|\mu|+|\nu|} \frac{|\{G_1 \subseteq
G_{\lambda}:G_{\lambda}/G_1=\mu,G_1=\nu\}|}{p^{\sum (\lambda_i')^2}
\prod_i (\frac{1}{p})_{m_i(\lambda)}} = \frac{1} {p^{\sum (\mu_i')^2}
\prod_i (\frac{1}{p})_{m_i(\mu)}} \frac{1}{p^{\sum (\nu_i')^2} \prod_i
(\frac{1}{p})_{m_i(\nu)}}. \] Fixing $\mu$, summing the left hand side over all $\nu$ of size
$n-|\mu|$, and applying Lemma \ref{duality} yields

\begin{eqnarray*}
\sum_{\lambda} \sum_{\nu} \frac{|\{G_1 \subseteq
G_{\lambda}:G_{\lambda}/G_1=\mu,G_1=\nu\}|}{p^{\sum (\lambda_i')^2}
\prod_i (\frac{1}{p})_{m_i(\lambda)}} & = & \sum_{\lambda} \sum_{\nu} \frac{|\{G_1 \subseteq
G_{\lambda}:G_{\lambda}/G_1=\nu,G_1=\mu\}|}{p^{\sum (\lambda_i')^2}
\prod_i (\frac{1}{p})_{m_i(\lambda)}}\\
& = & \sum_{\lambda} \frac{|\{G_1 \subseteq
G_{\lambda}:G_1=\mu\}|}{p^{\sum (\lambda_i')^2}
\prod_i (\frac{1}{p})_{m_i(\lambda)}}.
\end{eqnarray*} Fixing $\mu$, summing the right hand side over all $\nu$ of size
$n-|\mu|$, and applying Lemma \ref{prod} gives that

\[ \frac{1} {p^{\sum (\mu_i')^2} \prod_i (\frac{1}{p})_{m_i(\mu)}}
\sum_{\nu \vdash n-|\mu|} \frac{1}{p^{\sum (\nu_i')^2} \prod_i
(\frac{1}{p})_{m_i(\nu)}} =  \frac{1}{p^{\sum (\mu_i')^2} \prod_i
(\frac{1}{p})_{m_i(\mu)}} \frac{1}{p^{n-|\mu|}
(\frac{1}{p})_{n-|\mu|}},\] proving the lemma. \end{proof}

\subsection{Upper triangular matrices over a finite field} Let
$T(n,p)$ denote the set of upper triangular elements of $GL(n,p)$ with
$1$'s along the main diagonal. From the theory of wild quivers there
is a provable sense in which the conjugacy classes of $T(n,p)$ cannot
be classified. Nevertheless, as emerges from work of Kirillov
\cite{K1,K2} and Borodin \cite{B}, it is interesting to study the
Jordan form of elements of $T(n,p)$. As with the unipotent conjugacy
classes of $GL(n,p)$, the possible Jordan forms correspond to
partitions $\lambda$ of $n$.

	Theorem \ref{express} gives five expressions for the
probability that an element of $T(n,p)$ has Jordan form of type
$\lambda$. As is evident from the proof, most of the hard work at the
heart of these formulas has been carried out by others. Nevertheless,
at least one of these expressions is useful, and to the best of our
knowledge none of these formulas has appeared elsewhere. $P_{\lambda}$ will denote the Hall-Littlewood polynomial
of the previous subsection. By a standard Young tableau $S$ of size
$|S|=n$ is meant an assignment of $\{1,\cdots,n\}$ to the dots of the
partition such that each of $\{1,\cdots,n\}$ appears exactly once, and
the entries increase along the rows and columns. For instance, \[
\begin{array}{c c c c c}
                1 & 3 & 5 & 6 &   \\
                2 & 4 & 7 &  &    \\
                8 & 9 &  &  &    
          \end{array}  \] is a standard Young tableau.

\begin{theorem} \label{express} The probability that a uniformly
chosen element of $T(n,p)$ has Jordan form of type $\lambda$ is equal
to each of the following:

\begin{enumerate}
\item $\frac{(p-1)^n
P_{\lambda}(\frac{1}{p},\frac{1}{p^2},\frac{1}{p^3},\cdots;\frac{1}{p})
fix_{\lambda}(p)}{p^{n(\lambda)}}$, where $fix_{\lambda}(p)$ is the
number of complete flags of an $n$-dimensional vector space over a
field of size $p$ which are fixed by a unipotent element $u$ of type
$\lambda$.

\item $\frac{(p-1)^n P_{\lambda}(\frac{1}{p},
\frac{1}{p^2},\frac{1}{p^3},\cdots;\frac{1}{p})
Q_{(1)^n}^{\lambda}(p)}{p^{n(\lambda)}}$, where
$Q_{(1)^n}^{\lambda}(p)$ is a Green's polynomial as defined on page
247 of \cite{Mac}.

\item $(p-1)^n P_{\lambda}(\frac{1}{p} ,\frac{1}{p^2},
\frac{1}{p^3},\cdots;\frac{1}{p}) \sum_{\mu} dim(\chi^{\mu})
K_{\mu,\lambda}(\frac{1}{p})$, where $\mu$ is a partition of $n$,
$dim(\chi^{\mu})$ is the dimension of the irreducible representation of
$S_n$ of type $\mu$, and $K_{\mu,\lambda}$ is the Kostka-Foulkes
polynomial.

\item $\frac{(p-1)^n P_{\lambda}(\frac{1}{p},\frac{1}{p^2},\frac{1}{p^3},\cdots;\frac{1}{p})
chain_{\lambda}(p)}{p^{n(\lambda)}}$, where $chain_{\lambda}(p)$ is
the number of maximal length chains of subgroups in an abelian $p$-group of
type $\lambda$.

\item
$P_{\lambda}(1,\frac{1}{p},\frac{1}{p^2},\frac{1}{p^3},\cdots;\frac{1}{p})
\sum_{S} \prod_{j=1}^n (1-\frac{1}{p^{m^*(\Lambda_j)}})$, where
the sum is over all standard Young tableaux of shape $\lambda$, and
$m^*(\Lambda_j)$ is the number of parts in the subtableau formed
by $\{1,\cdots,j\}$ which are equal to the column number of $j$.
\end{enumerate}
\end{theorem}

\begin{proof} For the first assertion, observe that complete flags
correspond to cosets $GL(n,p)/B(n,p)$ where $B(n,p)$ is the subgroup
of all invertible upper triangular matrices. Note that $u \in GL(n,p)$
fixes the flag $gB(n,p)$ exactly when $g^{-1}ug \in B(n,p)$. The
unipotent elements of $B(n,p)$ are precisely $T(n,p)$. Thus the number
of complete flags fixed by $u$ is $\frac{1}{(p-1)^n |T(n,p)|} |\{g:
g^{-1}ug \in T(n,p)\}|$. It follows that the sought probability is
equal to $(p-1)^n fix_{\lambda}(p)$ multiplied by the probability that
an element of $GL(n,p)$ is unipotent of type $\lambda$. The first
assertion then follows from Lemma \ref{classsize}.

	The second part follows from the first part since by page 187
of \cite{Mac}, $Q_{(1)^n}^{\lambda}(p)$ is the number of complete
flags of an $n$-dimensional vector space over a field of size $p$
which are fixed by a unipotent element of type $\lambda$. The third
part follows from the second part and a formula for
$Q_{(1)^n}^{\lambda}(p)$ on page 247 of \cite{Mac}. The fourth part
follows from the third part and a formula for $\sum_{\mu}
dim(\chi^{\mu}) K_{\mu,\lambda}(\frac{1}{p})$ in \cite{Kir}.

	For the fifth assertion, a result on page 197 of \cite{Mac}
gives that the number of maximal length chains of subgroups in an
abelian $p$-group of type $\lambda$ is equal to $\frac{p^{n(\lambda)}}
{(1-\frac{1}{p})^n} \sum_{S} \prod_{j=1}^n
(1-\frac{1}{p^{m^*(\Lambda_j)}})$. Observing that for a partition
$\lambda$ of $n$, $P_{\lambda}(1,\frac{1}{p},
\frac{1}{p^2},\frac{1}{p^3},\cdots;\frac{1}{p})=p^n
P_{\lambda}(\frac{1}{p}, \frac{1}{p^2},\frac{1}{p^3},
\cdots;\frac{1}{p})$, the result follows.  \end{proof}

	As a corollary of Theorem \ref{express}, we recover the
``Division Algorithm'' of Borodin \cite{B}, which gives a
probabilistic way of growing partitions a dot at a time such
that the chance of getting $\lambda$ after $n$ steps is equal to the
chance that a uniformly chosen element of $T(n,p)$ has Jordan type
$\lambda$. We include our proof as it uses symmetric functions, which
aren't mentioned in the literature on probability in the upper
triangular matrices.

	We remark that a wonderful application of the division
algorithm was found by Borodin \cite{B}, who proved asymptotic
normality theorems for the lengths of the longest parts of the
partition corresponding to a random element of $T(n,p)$, and even
found the covariance matrix. We give another application in Section
\ref{applictriang}.
	
\begin{cor} (\cite{B}) Starting with the empty partition, at each
step transition from a partition $\lambda$ to a partition $\Lambda$ by
adding a dot to column $i$ chosen according to the rules

\begin{itemize}
\item $ i=1$ with  probability   $\frac{1}{p^{\lambda_1'}}$
\item $i=j>1$ with probability  $\frac{1}{p^{\lambda_j'}}-\frac{1}{p^{\lambda_{j-1}'}}$
\end{itemize}
\end{cor}

\begin{proof} For a standard Young tableau $S$, let $\Lambda_j(S)$ be
the partition formed by the entries $\{1,\cdots,j\}$ of $S$. It
suffices to prove that at step $j$ the division algorithm goes from
$\Lambda_{j-1}$ to $\Lambda_j$ with probability
$\frac{P_{\Lambda_j}(1,\frac{1}{p},\frac{1}{p^2},
\frac{1}{p^3},\cdots;\frac{1}{p})}{P_{\Lambda_{j-1}}(1,\frac{1}{p},\frac{1}{p^2},
\frac{1}{p^3} ,\cdots;\frac{1}{p})} (1-\frac{1}{p^{m^*(\Lambda_j)}})$,
because then the probability the Borodin's algorithm gives $\lambda$
at step $n=|\lambda|$ is

\[ \sum_{S : shape(S)=\lambda} \prod_{j=1}^n
\frac{P_{\Lambda_j}(1,\frac{1}{p},\frac{1}{p^2},
\frac{1}{p^3},\cdots;\frac{1}{p})}{P_{\Lambda_{j-1}}(1,\frac{1}{p},\frac{1}{p^2},
\frac{1}{p^3} ,\cdots;\frac{1}{p})} (1-\frac{1}{p^{m^*(\Lambda_j)}}) =
P_{\lambda}(1,\frac{1}{p},\frac{1}{p^2},\frac{1}{p^3},\cdots;\frac{1}{p})
\sum_{S} \prod_{j=1}^n (1-\frac{1}{p^{m^*(\Lambda_j)}}),\] as desired
from part 5 of Theorem \ref{express}. The fact that the division
algorithm
 goes from
$\Lambda_{j-1}$ to $\Lambda_j$ with probability
$\frac{P_{\Lambda_j}(1,\frac{1}{p},\frac{1}{p^2},
\frac{1}{p^3},\cdots;\frac{1}{p})}{P_{\Lambda_{j-1}}(1,\frac{1}{p},\frac{1}{p^2},
\frac{1}{p^3} ,\cdots;\frac{1}{p})} (1-\frac{1}{p^{m^*(\Lambda_j)}})$ follows, after algebraic manipulations, from Macdonald's principle specialization
formula (page 337 of \cite{Mac})

\[ P_{\lambda}(1,\frac{1}{p},\frac{1}{p^2},
\frac{1}{p^3},\cdots;\frac{1}{p}) = p^{n+n(\lambda)} \prod_i
\frac{1}{p^{\lambda_i'^2} (\frac{1}{p})_{m_i(\lambda)}}.\]
\end{proof}

As a remark, we observe that the division algorithm ties in
with an algorithm of the author for growing random parititions
distributed according to the $n \rightarrow \infty$ law of the
partition corresponding to the polynomial $z-1$ in the Jordan form of
a random element of $GL(n,p)$. One version of that algorithm \cite{F}
is

\begin{description}

\item [Step 0] Start with $\lambda$ the empty partition and
$N=1$. Also start with a collection of coins indexed by the natural
numbers such that coin $i$ has probability $\frac{1}{p^i}$ of heads and
probability $1-\frac{1}{p^i}$ of tails.

\item [Step 1] Flip coin $N$.

\item [Step 2a] If coin $N$ comes up tails, leave $\lambda$ unchanged,
set $N=N+1$ and go to Step 1.

\item [Step 2b] If coin $N$ comes up heads, let $j$ be the number of
the last column of $\lambda$ whose size was increased during a toss of
coin $N$ (on the first toss of coin $N$ which comes up heads, set
$j=0$). Pick an integer $S>j$ according to the rule that $S=j+1$ with
probability $\frac{1}{p^{\lambda_{j+1}'}}$ and $S=s>j+1$ with
probability $\frac{1}{p^{\lambda_s'}} - \frac{1}{p^{\lambda_{s-1}'}}$
otherwise. Then increase the size of column $S$ of $\lambda$ by 1 and
go to Step 1.

\end{description}

{\bf Remarks:}
\begin{enumerate}

\item Probabilistic algorithms similar to the one just described for
$GL(n,p)$ were used profitably in \cite{F} to prove group theoretic
results of Lusztig, Rudvalis/Shinoda, and Steinberg typically proved
by techniques such as character theory or Moebius inversion.

\item Observe that if one condition on the (probability zero) event
that each coin comes up heads exactly once, the transition rules are
the same as for the division lemma.
\end{enumerate}

\section{Applications} \label{applications} In this section we return
to the problem which motivated this paper: studying the
eigenvalues of unipotent matrices in the permutation representation on
lines. Lemma \ref{translate} describes the cycle structure of the
permutation action of a unipotent element $\alpha$ of $GL(n,p)$ on
lines in $V$ in terms of the partition parametrizing the conjugacy
class of $\alpha$.

\begin{lemma} \label{translate} Let $\alpha$ be a unipotent element of
$GL(n,p)$ with conjugacy class of type $\lambda$. Every orbit of the
action of $\alpha$ on the set of lines in $V$ has size $p^{r}$ for
some $r \geq 0$. The number of orbits of size $p^{r}$ is

\[\begin{array}{ll}
	 \frac{p^{\lambda_1'+\cdots+\lambda_{p^r}'}-
p^{\lambda_1'+\cdots+\lambda_{p^{r-1}}'}}{p-1}& \mbox{if $r \geq 1$}\\		\frac{p^{\lambda_1'}-1}{p-1} & \mbox{if $r=0$}.	\end{array} \]

\end{lemma}

\begin{proof}As discussed at the beginning of Section \ref{ident}, the
matrix $\alpha$ may be assumed to be

\[ \left( \begin{array}{c c c c}
		M_{\lambda_1} & 0 & 0 &0 \\
		0 & M_{\lambda_2} & 0 & 0\\
		0 & 0 & M_{\lambda_3} & \cdots \\
		0 & 0 & 0 & \cdots
	  \end{array} \right), \] where $M_i$ is the $i*i$ matrix of the form

\[ \left( \begin{array}{c c c c c c}
		1 & 1 & 0 & \cdots & \cdots & 0\\
		0 & 1 & 1 & 0 & \cdots & 0\\
		0 & 0 & 1 & 1 & \cdots & 0\\
		\cdots & \cdots & \cdots & \cdots & \cdots & \cdots \\
		\cdots & \cdots & \cdots & 0 & 1&1\\
		0 & 0 & 0 & \cdots & 0 & 1
	  \end{array} \right). \] Let $E_i=M_i-Id$, where $Id$ is the
identity matrix.

From this explicit form all eigenvalues of $\alpha^m, m \geq 0$ are
$1$. Thus if $\alpha^m$ fixes a line, it fixes it pointwise. Hence the
number of lines fixed by $\alpha^m$ is one less than the number of
points it fixes, all divided by $p-1$, and we are reduced to studying
the action of $\alpha$ of non-zero vectors. It is easily proved that $M_i$ has order $p^a$, where $p^{a-1} < i
\leq p^a$. Hence if $\alpha^m(x_1,\cdots,x_n)=(x_1,\cdots,x_n)$ with
some $x_i$ non-zero, and $m$ is the smallest non-negative integer with
this property, then $m$ is a power of $p$. Thus all orbits of $\alpha$
on the lines of $V$ have size $p^r$ for $r \geq 0$.

We next claim that $\alpha^{p^r}, r \geq 0$ fixes a vector

\[ (x_1,\cdots,x_{\lambda_1}, x_{\lambda_1+1}, \cdots,
x_{\lambda_1+\lambda_2},x_{\lambda_1+\lambda_2+1},\cdots,x_{\lambda_1+\lambda_2+\lambda_3},\cdots,x_n) \] if and only if

\[ x_{\lambda_1+\cdots+\lambda_{i-1}+p^a+1} =
x_{\lambda_1+\cdots+\lambda_{i-1}+p^a+2} = \cdots = x_{\lambda_1+\cdots+\lambda_i}=0 \ for \ i: \lambda_i > p^a.\] It suffices to prove this claim when $\lambda$ has one part
$\lambda_1$ of size $n$. Observe that the $i$th coordinate of
$\alpha^{p^a} (x_1,\cdots,x_n)$ is $\sum_{j=i}^n {p^a \choose j-i}
x_j$. If $n \leq p^a$, then $\alpha^{p^a}$ fixes all
$(x_1,\cdots,x_n)$ because ${p^a \choose r} = 0 \ mod \ p$ for
$r<p^a$. If $n>p^a$, then $\alpha^{p^a}$ fixes all $(x_1,\cdots,x_n)$
such that $x_{p^a+1}=\cdots=x_n=0$, for the same reason. Finally, if
$n>p^a$ and $x_j \neq 0$ for some $j>p^a$, let $j$ be the largest such
subscript. Then the $j-p^a$th coordinate of $\alpha^{p^a}
(x_1,\cdots,x_n)$ is equal to $x_{j-p^a} + x_j$ mod $p$, showing that
$\alpha^{p^a}$ does not fix such $(x_1,\cdots,x_n)$.

This explicit description of fixed vectors (hence of fixed lines) of
$\alpha^{p^a}$ yields the formula of the lemma for $r \geq 1$, because
the number of lines in an orbit of size $p^r$ is the difference
between the number of lines fixed by $\alpha^{p^r}$ and the number of
lines fixed by $\alpha^{p^{r-1}}$. The formula for the number of lines
in an orbit of size 1 follows because there are a total of
$\frac{p^n-1}{p-1}$ lines.\end{proof}

\subsection{Unipotent elements of $GL(n,p)$} \label{applicunip}

 Let $\alpha$ be a uniformly chosen unipotent element of
$GL(n,p)$. Each element of $GL(n,p)$ permutes the lines in $V$ and
thus defines a permutation matrix, which has complex eigenvalues. Each
size $p^r$ orbit of $\alpha$ on lines gives $p^r$ eigenvalues, with
one at each of the $p^r$th roots of unity. For $\theta \in (0,1)$, define a random variable
$X^{\theta}$ by letting $X^{\theta}(\alpha)$ be the number of
eigenvalues of $\alpha$ in the interval $(1,e^{2 \pi i \theta}]$ on
the unit circle. For $r \geq 1$, define random variables $X_r$ on the
unipotent elements of $GL(n,p)$ by

\[ X_r(\alpha) =
\frac{p^{\lambda_1'(\alpha)+\cdots+\lambda_r'(\alpha)}-p^{\lambda_1'(\alpha)+\cdots+\lambda_{r-1}'(\alpha)}}{p-1}. \] Clearly $X_r(\alpha)=0$ if $r>n$. Let $\lfloor y \rfloor$ denote the greatest integer less than
$y$. Lemma \ref{translate} implies that

\[ X^{\theta} = X_1 \lfloor \theta \rfloor + \sum_{r \geq 1}
\frac{X_{p^{r-1}+1} + \cdots + X_{p^r}}{p^r} \lfloor p^r \theta
\rfloor.\] This relationship (analogous to one used in \cite{Wieand})
will reduce the computation of the mean of $X^{\theta}$ to similar
computations for the random variables $X_r$, which will now be carried
out.

Let $E_n$ denote the expected value with respect to the uniform
distribution on the unipotent elements of $GL(n,p)$.

\begin{theorem} \label{mean} For $1 \leq r \leq n$, \[ E_n(X_r) = \frac{p^r
(1-\frac{1}{p^{n-r+1}})\cdots(1-\frac{1}{p^n})}{p-1} .\] \end{theorem}

\begin{proof} By Lemma \ref{formula},

\[ E_n(X_r) = \sum_{\lambda \vdash n} \frac{p^n
(\frac{1}{p})_n}{p^{\sum (\lambda_i')^2} \prod_i
(\frac{1}{p})_{m_i(\lambda)}}
\frac{p^{\lambda_1'(\alpha)+\cdots+\lambda_r'(\alpha)}
-p^{\lambda_1'(\alpha)+\cdots+\lambda_{r-1}'(\alpha)}}{p-1}.\] Observe that $\frac{p^{\lambda_1'+\cdots+\lambda_r'}-
p^{\lambda_1'+\cdots+\lambda_{r-1}'}}{p^r-p^{r-1}}$ is the number of
subgroups of $G_{\lambda}$ of type $\nu=(r)$. This is because the
total number of elements of order $p^r$ in $G_{\lambda}$ is
$p^{\lambda_1'+\cdots+\lambda_r'}-p^{\lambda_1'+\cdots+\lambda_{r-1}'}$,
and every subgroup of type $\nu=(r)$ has $p^r-p^{r-1}$
generators. Therefore, using Lemma \ref{likeMac},

\begin{eqnarray*}
E_n(X_r) &=& p^n (\frac{1}{p})_n \frac{p^r-p^{r-1}}{p-1} \sum_{\lambda \vdash n} \frac{|\{G_1 \subseteq
G_{\lambda} : G_1=(r)\}|}{p^{\sum (\lambda_i')^2} \prod_i (\frac{1}{p})_{m_i(\lambda)}}\\
& = & \left(p^n (\frac{1}{p})_n \frac{p^r-p^{r-1}}{p-1} \right) \left(\frac{1}{p^r(1-\frac{1}{p})} \frac{1}{p^{n-r} (\frac{1}{p})_{n-r}} \right)\\
& = & \frac{p^r (1-\frac{1}{p^{n-r+1}})\cdots(1-\frac{1}{p^n})}{p-1}.
\end{eqnarray*}
\end{proof}  

	Corollary \ref{mean2} uses Theorem \ref{mean} to bound the
mean of $X^{\theta}$.

\begin{cor} \label{mean2} $E_n(X^{\theta}) = \theta
\frac{p^n-1}{p-1} - O(\frac{p^n}{n})$.
\end{cor}

\begin{proof} Let $\{y\}=y-\lfloor y \rfloor$ denote the fractional
part of a positive number $y$. Theorem \ref{mean} and the writing of
$X^{\theta}$ in terms of the $X_r$'s imply that

\begin{eqnarray*}
E_n(X^{\theta}) & = & \theta E_n(\sum_{i \geq 1} X_i)
- \sum_{r \geq 1} \{p^r \theta\} E_n
(\frac{X_{p^{r-1}+1}+\cdots+X_{p^r}}{p^r})\\
& = & \theta
\frac{p^n-1}{p-1} - \sum_{r \geq 1}
\{p^r \theta\} E_n (\frac{X_{p^{r-1}+1}+\cdots+X_{p^r}}{p^r})\\
& \geq
& \theta \frac{p^n-1}{p-1} - \sum_{r \geq 1} E_n
(\frac{X_{p^{r-1}+1}+\cdots+X_{p^r}}{p^r})\\
& \geq & \theta \frac{p^n-1}{p-1} - (\sum_{r=1}^{\lfloor log_p(n) \rfloor)}
\frac{p^{p^{r-1}+1} + \cdots + p^{p^r}}{(p-1)p^r}) -
(\frac{p^{p^{\lfloor log_p(n) \rfloor}+1} + \cdots +
p^n}{(p-1)p^{\lfloor log_p(n) \rfloor +1}}).
\end{eqnarray*}

We suppose for simplicity that $n \neq p^{p^r}+1$ for some $r$ (the
case $n = p^{p^r}+1$ is similar). Continuing,

\begin{eqnarray*}
E_n(X^{\theta}) & \geq &\theta \frac{p^n-1}{p-1} - (\sum_{r=1}^{\lfloor log_p(n) \rfloor} \frac{p^{p^r+1}}{(p-1)^2 \frac{n}{p^{\lfloor log_p(n) \rfloor -r +1}}}) - \frac{p^{n+1}}{(p-1)^2 n} =  \theta \frac{p^n-1}{p-1} - O(\frac{p^n}{n}).
\end{eqnarray*}
\end{proof}

	The approach here appears to extend to the computation of
higher moments, but the computations are formidable. For example one
can show that if $1 \leq r \leq s \leq n$, then \[ E_n(X_rX_s) =
\frac{p^{r+s-1}}{p-1} [\frac{p}{p-1}(1-\frac{1}{p^{n-s-r+1}}) \cdots
(1-\frac{1}{p^n}) + \sum_{a=0}^{r-1} (1-\frac{1}{p^{n-a-s+1}}) \cdots
(1-\frac{1}{p^n})] . \]

\subsection{Upper triangular matrices over a finite field}
\label{applictriang}

 Let $\alpha$ be a uniformly chosen element of $T(n,p)$. Recall that
$\alpha$ is unipotent by the definition of $T(n,p)$. Each element of
$T(n,p)$ permutes the lines in $V$ and thus defines a permutation
matrix, which has complex eigenvalues. Each size $p^r$ orbit of
$\alpha$ on lines gives $p^r$ eigenvalues, with one at each of the
$p^r$th roots of unity. For $\theta \in (0,1)$, define a random
variable $X^{\theta}$ by letting $X^{\theta}(\alpha)$ be the number of
eigenvalues of $\alpha$ in the interval $(1,e^{2 \pi i \theta}]$ on
the unit circle. For $r \geq 1$, define random variables $X_r$ on the
unipotent elements of $T(n,p)$ by

\[ X_r(\alpha) =
\frac{p^{\lambda_1'(\alpha)+\cdots+\lambda_r'(\alpha)}-p^{\lambda_1'(\alpha)+\cdots+\lambda_{r-1}'(\alpha)}}{p-1}. \] Let $\lfloor y \rfloor$ denote the greatest integer less than
$y$. Lemma \ref{translate} implies that

\[ X^{\theta} = X_1 \lfloor \theta \rfloor + \sum_{r \geq 1}
\frac{X_{p^{r-1}+1} + \cdots + X_{p^r}}{p^r} \lfloor p^r \theta
\rfloor.\] As for the case of $GL(n,p)$ this relationship will reduces
the computation of the mean of $X^{\theta}$ to similar computations
for the random variables $X_r$.

Let $E_n$ denote the expected value with respect to the uniform
distribution on the unipotent elements of $T(n,p)$. Theorem \ref{mean3}
shows that the expected value of $X_r$ is surprisingly simple. As one
sees from the case $p=2$, the result is quite different from that
of Theorem \ref{mean}.

\begin{theorem} \label{mean3} For $1 \leq r \leq n$, \[ E_n(X_r) =
(p-1)^{r-1} {n \choose r} .\] \end{theorem}

\begin{proof} We proceed by joint induction on $n$ and $r$, the base
case $n=r=1$ being clear. Let $Prob(S)$ denote the probability that
Borodin's growth algorithm yields the standard Young tableau $S$
after $|S|$ steps. Let $col(n)$ be the column number of $n$ in
$S$. With all sums being over standard Young tableaux, observe that

\begin{eqnarray*}
E_n(p^{\lambda_1'+\cdots+\lambda_r'}) & = & \sum_{S: |S|=n} p^{\lambda_1'(S)+\cdots+\lambda_r'(S)} Prob(S)\\
& = & \sum_{S: |S|=n,col(n)=1} p^{\lambda_1'(S)+\cdots+\lambda_r'(S)} Prob(S)\\&& + \sum_{S:|S|=n,1<col(n)=j\leq r} p^{\lambda_1'(S)+\cdots+\lambda_r'(S)} Prob(S)\\ &&+ \sum_{S:|S|=n,col(n)>r} p^{\lambda_1'(S)+\cdots+\lambda_r'(S)} Prob(S)\\
& = & \sum_{S':|S'|=n-1} p^{\lambda_1'(S')+\cdots+\lambda_r'(S')+1} Prob(S') \frac{1}{p^{\lambda_1'(S')}}\\&& + \sum_{j=2}^r \sum_{S':|S'|=n-1} p^{\lambda_1'(S')+\cdots+\lambda_r'(S')+1} Prob(S') (\frac{1}{p^{\lambda_j'(S')}}-\frac{1}{p^{\lambda_{j-1}'(S')}})\\&&+ \sum_{j>r}  \sum_{S':|S'|=n-1} p^{\lambda_1'(S')+\cdots+\lambda_r'(S')} Prob(S') (\frac{1}{p^{\lambda_j'(S')}}-\frac{1}{p^{\lambda_{j-1}'(S')}})\\
& = & p E_{n-1}(p^{\lambda_2'+\cdots+\lambda_r'}) + p E_{n-1}(p^{\lambda_1'+\cdots+\lambda_{r-1}'}-p^{\lambda_2'+\cdots+\lambda_r'})\\&& + E_{n-1}(p^{\lambda_1'+\cdots+\lambda_{r}'}-p^{\lambda_1'+\cdots+\lambda_{r-1}'})\\
& = & (p-1) E_{n-1}(p^{\lambda_1'+\cdots+\lambda_{r-1}'}) + E_{n-1}(p^{\lambda_1'+\cdots+\lambda_{r}'})\\
& = & (p-1)^{r-1} {n-1 \choose r-1} + (p-1)^{r-1} {n-1 \choose r}\\
& = & (p-1)^{r-1} {n \choose r}.
\end{eqnarray*}
\end{proof}

	Corollary \ref{mean4} follows by using Theorem \ref{mean3} and
arguing along the lines of Corollary \ref{mean2}.

\begin{cor} \label{mean4} $\theta \frac{p^n-1}{p-1}-p \sum_{r=1}^n
\frac{(p-1)^{r-1} {n \choose r}}{r} \leq E_n(X^{\theta}) \leq \theta
\frac{p^n-1}{p-1}$. \end{cor}

\section{Acknowledgments} The author thanks Persi Diaconis for helpful
references. This research was supported by an NSF Postdoctoral
Fellowship.

\end{document}